\documentclass[12pt,a4paper,draft]{amsart}

\def\BB{\mathbb B}
\def\CC{\mathbb C}
\def\CO{\mathcal O}
\def\DD{\mathbb D}
\def\EE{\mathbb E}
\def\NN{\mathbb N}
\def\PSH{\mathcal{PSH}}
\def\RR{\mathbb R}
\def\too{\longrightarrow}

\theoremstyle{definition}
\newtheorem{definition}{Definition}[section]
\newtheorem{remark}[definition]{Remark}

\theoremstyle{plain}
\newtheorem{corollary}[definition]{Corollary}
\newtheorem{lemma}[definition]{Lemma}
\newtheorem{proposition}[definition]{Proposition}
\newtheorem{theorem}[definition]{Theorem}

\begin{document}

\title[M\"obius function in complex ellipsoids]{M\"obius function of coordinate hyperplanes in complex ellipsoids}
\author{Witold Jarnicki}
\address{Jagiellonian University, Institute of Mathematics\newline
\indent Reymonta 4, 30-059 Krak\'ow, Poland}
\email{wmj@im.uj.edu.pl}
\curraddr{Universit\"at Osnabr\"uck, Fachbereich Mathematik/Informatik\newline
\indent Albrechtstra\ss e 28, 49069 Osnabr\"uck, Germany}
\thanks{The author was supported in part by KBN grant no.~2 P03A 015 22.}

\begin{abstract}
For $p_1,\dots,p_n>0$, let $\EE=\{z\in\CC^n:\sum_{j=1}^n|z_j|^{2p_j}<1\}$ be
a complex ellipsoid. We present effective formulas for the generalized
M\"obius and Green functions $m_{\EE}(A,\cdot)$, $g_{\EE}(A,\cdot)$ in the case
where $A:=\{z\in\EE: z_1\cdot\dots\cdot z_k=0\}$ ($1\leq k\leq n$).
%\footnote{}
\end{abstract}

\maketitle

\section{Introduction and main results}

Let $G\subset\CC^n$ be a domain and let $A\subset G$. Define the 
{\it generalized Green function}
\begin{multline*}
g_G(A,z):=\sup\{u(z):\quad u:G\too[0,1),\;\log u\in\PSH(G),\\
\forall_{a\in A}\;\exists_{C=C(u,a)>0}\;\forall_{w\in G}:\;u(w)\leq C\|w-a\|\},\quad
z\in G,
\end{multline*}
and the {\it generalized M\"obius function}
$$
m_G(A,z):=\sup\{|f(z)|:f\in\CO(G,\DD),\;f|_A\equiv0\},\quad z\in G,
$$
where $\DD$ denotes the unit disc (cf. \cite{JarJarPfl2003}).

The generalized M\"obius and Green functions were recently studied by many
authors (e.g. \cite{LarSig1998}). Despite various properties proven, 
effective formulas for
$m_G(A,\cdot)$ and $g_G(A,\cdot)$ are known only in a few special cases.
An effective formula for the generalized Green function in the case of
Euclidean ball $G=\BB_2\subset\CC^2$ and union of coordinate lines
$A=\{(z,w)\in\BB_2: zw=0\}$
was given in \cite{Ngu2003}.

Let $n\geq2$. Let $p_1,\dots,p_n>0$. Define
the {\it complex ellipsoid with exponents $p_1,\dots,p_n$}:
$$
\EE=\EE_{p_1,\dots,p_n}:=\Big\{z\in\CC^n:\sum_{j=1}^n|z_j|^{2p_j}<1\Big\}.
$$
Observe that $\EE_{p_1,\dots,p_n}$ is convex if and only if $p_j\geq1/2$,
$j=1,\dots,n$ (cf. \cite{JarPfl1993}, \S\;8.4).

For $G\subset\CC^n$, $k\in\{1,\dots,n\}$, we consider the sets
$$
A=A_{G,k}:=\{z\in G:z_1\cdot\dots\cdot z_k=0\}.
$$

We give an effective formula for $g_{\EE}(A_{\EE,k},\cdot)$ for any complex
ellipsoid $\EE\subset\CC^n$ and any $k\in\{1,\dots,n\}$ (Theorem~\ref{one}(a)).
This is a generalization of a result in \cite{Ngu2003}
(with a much simpler proof). We prove that the
$m_{\EE}(A_{\EE,k},\cdot)\equiv g_{\EE}(A_{\EE,k},\cdot)$ for the cases where
$\EE$ is convex (Theorem~\ref{one}(b)) or where $k=1$, $n=2$, $p_2\geq1/2$ (Theorem~\ref{one}(c)).
This is not true in general nonconvex case; we prove that
$m_\EE(A_{\EE,k},\cdot)\not\equiv g_{\EE}(A_{\EE,k},\cdot)$ if there exists a $j\in\{k+1,\dots,n\}$ such
that $p_j<1/2$ (Theorem~\ref{one}(d)). Partial results for the case $k=n=2$ are given
(Theorem~\ref{one}(e)).

\begin{definition}\label{R}
Let $p_1,\dots,p_n>0$, $1\leq k\leq n$. Put
$\EE:=\EE_{p_1,\dots,p_n}$, $A:=A_{\EE,k}$.
Let $z\in\EE$ be such that the sequence
$(p_j|z_j|^{2p_j})_{j=1}^k$ is monotonically increasing. For $s\in\{1,\dots,k\}$ define
$$
q_s:=\sum_{j=1}^s(2p_j)^{-1},\quad
r_s(z):=1-\sum_{j=s+1}^n|z_j|^{2p_j},\quad
c_s(z):=r_s(z)/q_s.
$$
Let $d:=\max\{s\in\{1,\dots,k\}:2p_s|z_s|^{2p_s}\leq c_s(z)\}$. Define
\begin{equation*}
R_\EE(A,z)=
\prod_{j=1}^d|z_j|\Big(\frac{2p_j}{c_d(z)}\Big)^{\!\frac1{2p_j}}.\tag{$\ast$}
\end{equation*}
\end{definition}

The main results of the paper are the following:

\begin{theorem}\label{one}Under the above assumptions we have:
\begin{enumerate}
\item[(a)]$g_\EE(A,\cdot)\equiv R_\EE(A,\cdot)$,
\item[(b)]$m_\EE(A,\cdot)\equiv g_\EE(A,\cdot)\equiv R_\EE(A,\cdot)$, for $p_j\geq1/2$,
$j=1,\dots,n$,
\item[(c)]$m_\EE(A,\cdot)\equiv g_\EE(A,\cdot)\equiv R_\EE(A,\cdot)$, for $k=1$, $n=2$,
$p_2\geq1/2$,
\item[(d)]$m_\EE(A,\cdot)\not\equiv g_\EE(A,\cdot)$, if there exists a
$j\in\{k+1,\dots,n\}$ with $p_j<1/2$.
\item[(e)]$m_\EE(A,\cdot)\equiv g_\EE(A,\cdot)\equiv R_\EE(A,\cdot)$, for
$k=n=2$, $p_1\leq p_2$, and either $p_2\geq1/2$ or
$8p_1+4p_2(1-p_2)>1$.
\end{enumerate}
\end{theorem}

Observe that one can obtain the monotonicity of $(p_j|z_j|^{2p_j})_{j=1}^k$
for arbitrary $z\in\EE$ by permutation of coordinates. The final result is
a subdivision of $\EE$ into $2^k-1$ subsets, each with
$R_\EE(A,\cdot)$ given by a formula of type ($\ast$).

\begin{remark}
In the case where $p_1=\dots=p_n=1$, the domain $\EE$ is the
Euclidean ball $\BB_n$. Theorem \ref{one} may then be formulated
in the following simpler way.
\end{remark}

\begin{corollary}Assume that $|z_1|\leq\dots\leq|z_k|$ and let
$$
d:=\max\Big\{s\in\{1,\dots,k\}:s|z_s|^2+\sum_{j=s+1}^n|z_j|^2\leq1\Big\}.
$$
Then
$$
m_{\BB_n}(A,z)=g_{\BB_n}(A,z)=\Big(\frac d{1-\sum_{j=d+1}^n|z_j|^2}\Big)^{\!\frac d2}
\prod_{j=1}^d|z_j|.
$$
\end{corollary}

\begin{remark}It is unknown to the author whether
$m_\EE(A,\cdot)\equiv g_\EE(A,\cdot)$ for any $p_j>0$, $j=1,\dots,k$,
$p_j\geq1/2$,
$j=k+1,\dots,n$. However, Theorem~\ref{one}(e) and author's research
(cf. Remark~\ref{partial})
indicate that this is true at least in the case
$k=n=2$.
\end{remark}

\begin{remark}
Take $p_{j,\ell}\nearrow+\infty$, $j=1,\dots,n$. Then
$$
G_\ell:=\EE_{p_{1,\ell},\dots,p_{n,\ell}}\nearrow \DD^n.
$$
Consequently,
$$
m_{G_\ell}(A_{G_\ell,k},\cdot)=g_{G_\ell}(A_{G_\ell,k},\cdot)\searrow
m_{\DD^n}(A_{\DD^n,k},\cdot)=g_{\DD^n}(A_{\DD^n,k},\cdot)
$$
(cf. \cite{JarJarPfl2003},
Property 2.7) and hence
$$
m_{\DD^n}(A_{\DD^n,k},z)=|z_1|\cdot\dots\cdot|z_k|,\quad z\in\DD^n.
$$
However, we need the formula for $g_{\DD^n}(A_{\DD^n,k},\cdot)$ before
we prove Theorem~\ref{one} (cf. Lemma \ref{schwarz}).
\end{remark}

\section{Proof of Theorem~\ref{one}}

\begin{lemma}\label{schwarz}
Let $n\in\NN$, $1\leq k\leq n$. Then
$$
m_{\DD^n}(A_{\DD^n,k},z)=g_{\DD^n}(A_{\DD^n,k},z)=|z_1|\cdot\dots\cdot|z_k|,\quad
z\in\DD^n.
$$
\end{lemma}

\begin{proof}[Proof of Lemma~\ref{schwarz}] (Due to P.Pflug.) obviously,
$$
|z_1|\cdot\dots\cdot|z_k|\leq m_{\DD^n}(A_{\DD^n,k},z)\leq g_{\DD^n}(A_{\DD^n,k},z).
$$
It remains to prove that
$g_{\DD^n}(A_{\DD^n,k},z)\leq|z_1|\cdot\dots\cdot|z_k|$.

It is enough to prove that for any $1\leq k\leq n$ and for any function
$u:\DD^n\too[0,1)$ such
that $\log u\in\PSH(\DD^n)$ and $u(z)\leq C(a)\|z-a\|$, $a\in A_{\DD^n,k}$,
$z\in\DD^n$, we have
$u(z)\leq|z_1|\cdot\dots\cdot|z_k|$. We proceed by induction on $k$.

For $k=1$ the inequality follows from the Schwarz type lemma for
logarithmically subharmonic functions $u(\cdot,z_2,\dots,z_n)$,
$z_2,\dots,z_n\in \DD$.

For $k>1$ we first apply the case $k=1$ and get $u(z_1,\dots,z_n)\leq|z_1|$,
$z\in \DD^n$.
Applying the inductive assumption to $u(z_1,\cdot)/|z_1|$, $z_1\in \DD$, finishes
the proof.
\end{proof}

Take a $u:\EE\too[0,1)$ with $\log u\in\PSH(\EE)$ and
$u(\zeta)\leq C(a)\|\zeta-a\|$, $a\in A$, $\zeta\in\EE$.
Consider the mapping
\begin{multline*}
\iota_z:\DD^d\ni(\zeta_1,\dots,\zeta_d)\\\longmapsto
\Big(\zeta_1\Big(\frac{c_d(z)}{2p_1}\Big)^{\!\frac1{2p_1}},\dots,
\zeta_d\Big(\frac{c_d(z)}{2p_d}\Big)^{\!\frac1{2p_d}},
z_{d+1},\dots,z_n\Big)\in\EE.
\end{multline*}
Applying the holomorphic contractivity of the generalized Green function and
Lemma~\ref{schwarz} proves that $m_\EE(A,z)\leq g_\EE(A,z)\leq R_\EE(A,z)$.

It remains to prove that $m_\EE(A,z)\geq R_\EE(A,z)$ (resp. $g_\EE(A,z)\geq R_\EE(A,z)$).
In the case $d=k=n$ it suffices to take
\begin{align*}
f(\zeta):=\prod_{j=1}^d\zeta_j\Big(\frac{2p_j}{c_d(z)}\Big)^{\!\frac1{2p_j}}.
\end{align*}
Then $f\in\CO(\EE,\DD)$ and $|f(z)|=R_\EE(A,z)$.

It remains to prove Theorem~\ref{one} in the remaining cases. We may assume that $z_j\neq0$, $j=1,\dots,d$.
Consider the following lemma.

\begin{lemma}\label{max}Define $\EE':=\EE_{p_{d+1},\dots,p_n}$.
\begin{enumerate}
\item[(a)]
Let $v:\EE'\too[0,1)$, $v\not\equiv0$ be such
that $\log v\in\PSH(\EE')$ and $v(\zeta)\leq|\zeta_j|$,
$\zeta\in\EE'$, $j=d+1,\dots,k$, and the
mapping
$$
\EE'\ni(\zeta_{d+1},\dots,\zeta_n)\longmapsto v(\zeta_{d+1},\dots,\zeta_n)r_d^{q_d}(\zeta)
$$
attains its maximum value $M$ at $(z_{d+1},\dots,z_n)$.
Then $g_\EE(A,z)=R_\EE(A,z)$.

\item[(b)]The following conditions are equivalent:
\begin{itemize}
\item There exists an $h\in\CO(\EE')$, $h\not\equiv0$ such that $h(\zeta)=0$ for
$\zeta_{d+1}\cdot\dots\cdot\zeta_k=0$ and the mapping
$$
\EE'\ni(\zeta_{d+1},\dots,\zeta_n)\longmapsto|h(\zeta_{d+1},\dots,\zeta_n)|r_d^{q_d}(\zeta)
$$
attains its maximum value $M$ at the point $(z_{d+1},\dots,z_n)$,
\item $m_\EE(A,z)=R_\EE(A,z)$.
\end{itemize}
\end{enumerate}
\end{lemma}

We present a proof in Section~\ref{maxproof}. The result above reduces the proof of Theorem~\ref{one} to the following
propositions.

\begin{proposition}\label{a}There exists a $v$ as required in Lemma~\ref{max}(a).
\end{proposition}

\begin{proposition}\label{b}Using notation of Theorem~\ref{one}, assume that $p_j\geq1/2$, $j=d+1,\dots,n$.
Then there exists an $h\in\CO(\EE')$ as required in Lemma~\ref{max}(b).
\end{proposition}

\begin{proposition}\label{c}Using notation of Theorem~\ref{one}, assume that $p_{k+1}<1/2$. Then one cannot find an $h$ as
required in Lemma~\ref{max}(b) for $|z_\ell|\neq0$ small enough, $\ell=1,\dots,k+1$ and $z_\ell=0$, $\ell=k+2,\dots,n$.
\end{proposition}

\begin{proposition}\label{d}Using notation of Theorem~\ref{one}, assume that $n=k=2$. Additionally, assume that
$p_2\geq1/2$ or $8p_1+4p_2(1-p_2)>1$.
Then there exists an $h\in\CO(\DD)$ as required in Lemma~\ref{max}(b).
\end{proposition}

\section{Proof of Lemma~\ref{max}}\label{maxproof}

(a) Put
$$
u(\zeta_1,\dots,\zeta_n):=M^{-1}
\Big(\prod_{j=1}^d(2p_j)^{\frac1{2p_j}}|\zeta_j|\Big)q_d^{q_d}
v(\zeta_{d+1},\dots,\zeta_n).
$$
Obviously $\log u\in\PSH(\EE)$ and
$u(\zeta)\leq C|\zeta_j|\leq C\|\zeta-a\|$, $\zeta\in\EE$, whenever
$a_j=0$ for some $j\in\{1,\dots,k\}$. For $\zeta\in\EE$ we have:
$$
u(\zeta)\leq M^{-1}\Big(\frac{\sum_{j=1}^d|\zeta_j|^{2p_j}}{q_d}\Big)^{q_d}
q_d^{q_d}
v(\zeta_{d+1},\dots,\zeta_n)<1.\quad\footnotemark
$$
\footnotetext{Let $a_1,\dots,a_d\geq0$,
$w_1,\dots,w_d>0$. Then
$$
\prod_{j=1}^da_j^{w_j}\leq\left(\frac{\sum_{j=1}^dw_ja_j}{\sum_{j=1}^dw_j}\right)^{\sum_{j=1}^dw_j}.
$$}
Consequently $u:\EE\too[0,1)$. On the other hand:
$$
u(z)=M^{-1}R_\EE(A,z)r_d^{q_d}(z)v(z_{d+1},\dots,z_n)=R_\EE(A,z).
$$

(b) Assume that such an $h$ exists. Put
$$
f(\zeta_1,\dots,\zeta_n):=M^{-1}
\Big(\prod_{j=1}^d(2p_j)^{\frac1{2p_j}}\zeta_j\Big)q_d^{q_d}
h(\zeta_{d+1},\dots,\zeta_n).
$$
Observe that $f(\zeta)=0$ for $\zeta\in A$. Similarly as in (a) we prove that $|f(\zeta)|<1$,
$\zeta\in\EE$ and $f(z)=R_\EE(A,z)$.

Assume now that $m_\EE(A,z)=R_\EE(A,z)$. Let $f\in\CO(\EE,\DD)$ be such that $f|_A\equiv0$ and
$|f(z)|=R_\EE(A,z)$ (cf. \cite{JarJarPfl2003}, Property~2.5). Put
$$
h(\zeta):=\frac{\partial^df}{\partial z_1\dots\partial z_d}(0,\zeta).
$$
By definition, we have $h(\zeta)=0$ for $\zeta_{d+1}\cdot\dots\cdot\zeta_k=0$. Applying the Schwarz lemma
to the mapping $f\circ\iota_\zeta$, $\zeta\in\EE'$, we get
\begin{align*}
|h(\zeta_{d+1},\dots,\zeta_n)|r_d^{q_d}(\zeta)\leq&\Big(\prod_{j=1}^d(2p_j)^{\frac1{2p_j}}\Big)q_d^{q_d},
\quad\zeta\in\EE',\\
|h(z_{d+1},\dots,z_n)|r_d^{q_d}(z)=&\Big(\prod_{j=1}^d(2p_j)^{\frac1{2p_j}}\Big)q_d^{q_d}.
\end{align*}

\section{Proof of Propositions~\ref{a}, \ref{b}, \ref{c}, and \ref{d}}

\begin{proof}[Proof of Proposition~\ref{a}]We may assume that $z_{d+1},\dots,z_n\geq0$.
Consider functions of the form
$$
v_\alpha(\zeta_{d+1},\dots,\zeta_n)=
\Big(\prod_{j=d+1}^k|\zeta_j|^{1+\alpha_j}\Big)
\Big(\prod_{j=k+1}^n|\zeta_j|^{\alpha_j}\Big),
$$
where $\alpha=(\alpha_{d+1},\dots,\alpha_n)$, $\alpha_{d+1},\dots,\alpha_n\geq0$. Obviously $v:\EE'\too[0,1)$,
$\log v\in\PSH(\EE')$, and $v(\zeta)\leq|\zeta_j|\leq\|\zeta-a\|$,
$\zeta\in\EE'$, whenever $a_j=0$ for some $j\in\{d+1,\dots,k\}$. Since
$v_\alpha(\zeta_{d+1},\dots,\zeta_n)=
v_\alpha(|\zeta_{d+1}|,\dots,|\zeta_n|)$, it is enough to find an $\alpha$ such
that the function
$$
\EE'\cap\RR_+^{n-d}\ni(t_{d+1},\dots,t_n)\longmapsto
v_\alpha(t_{d+1},\dots,t_n)r_d^{q_d}(t)
$$
attains its maximum at $(z_{d+1},\dots,z_n)$. Considering the partial
(logarithmic) derivatives results in the following equations:
\begin{align*}
0=1+&\alpha_j-2p_jq_d\frac{t_j^{2p_j}}{r_d(t)},\quad j=d+1,\dots,k,\\
0=&\alpha_j-2p_jq_d\frac{t_j^{2p_j}}{r_d(t)},\quad j=k+1,\dots,n.
\end{align*}
These give formulas for $\alpha_{d+1},\dots,\alpha_n$ such that
$(z_{d+1},\dots,z_n)$ is the common zero of the derivatives. To prove that
there are no other points like this, consider a reformulation of the above
equations:
\begin{align*}
r_d(t)&=\frac{2p_jq_dt_j^{2p_j}}{1+\alpha_j},\quad j=d+1,\dots,k,\\
r_d(t)&=\frac{2p_jq_dt_j^{2p_j}}{\alpha_j},\quad j=k+1,\dots,n.
\end{align*}
The left side is decreasing in any of the variables, while the right sides
are increasing. Thus, at most one common zero is allowed.

It remains to check, whether $\alpha_j\geq0$, $j=d+1,\dots,n$. Obviously,
this is true for $j=k+1,\dots,n$ and in the remaining cases we have:
$$
\alpha_j=\frac{2p_jq_dt_j^{2p_j}-r_d(t)}{t_jr_d(t)}\geq0,
$$
since this is the way we have chosen $d$.
\end{proof}

\begin{proof}[Proof of Proposition~\ref{b}]
We may assume that $z_{d+1},\dots,z_n\geq0$.
Consider functions of the form
$$
h_\alpha(\zeta_{d+1},\dots,\zeta_n)=
\Big(\prod_{j=d+1}^k\zeta_je^{\alpha_j\zeta_j}\Big)
\Big(\prod_{j=k+1}^ne^{\alpha_j\zeta_j}\Big),
$$
where $\alpha=(\alpha_{d+1},\dots,\alpha_n)$,
$\alpha_{d+1},\dots,\alpha_n\geq0$. Since
$|h_\alpha(\zeta_{d+1},\dots,\zeta_n)|\leq
h_\alpha(|\zeta_{d+1}|,\dots,|\zeta_n|)$, it is enough to find an
$\alpha$ such that
$$
\EE'\cap\RR_+^{n-d}\ni(t_{d+1},\dots,t_n)\longmapsto
h_\alpha(t_{d+1},\dots,t_n)r_d^{q_d}(t)
$$
attains its maximum at $(z_{d+1},\dots,z_n)$. Considering the partial
(logarithmic) derivatives results in the following equations:
\begin{align*}
0=\frac1t_j+&\alpha_j-2p_jq_d\frac{t_j^{2p_j-1}}{r_d(t)},\quad j=d+1,\dots,k,\\
0=&\alpha_j-2p_jq_d\frac{t_j^{2p_j-1}}{r_d(t)},\quad j=k+1,\dots,n.
\end{align*}
We continue as in the proof of Proposition~\ref{a}.
\end{proof}

\begin{proof}[Proof of Proposition~\ref{c}]Consider the following two lemmas.

\begin{lemma}\label{neg1}Assume that there exist $0<c<b<1$ such that the
function
$$
\varphi:[0,b]\ni t\longmapsto(1-t^p)^{-q}
$$
is strictly concave and
$$
\varphi(0)+\frac bc(\varphi(c)-\varphi(0))>\varphi(b)+2.
$$
Let $f\in\CO(\DD)$, $f\not\equiv0$ such
that $|f(\zeta)|/\varphi(|\zeta|)$ attains its maximum at $w_0$. Then $w_0=0$
or $|w_0|\geq c$.
\end{lemma}

\begin{proof}Assume that $|w_0|\in(0,c)$. We may assume that
$|f(w_0)|=\varphi(|w_0|)$. Consider the function
$$
\psi:[0,b]\ni t\longmapsto|f(0)|+\frac t{|w_0|}|f(w_0)-f(0)|.
$$
From $\psi(0)\leq\varphi(0)$, $\psi(|w_0|)\geq\varphi(|w_0|)$, and the convexity condition we get:
\begin{align*}
\psi(b)=&|f(0)|+\frac b{|w_0|}|f(w_0)-f(0)|\\
    \geq&\varphi(0)+\frac b{|w_0|}|\varphi(|w_0|)-\varphi(0)|\\
    \geq&\varphi(0)+\frac bc|\varphi(c)-\varphi(0)|>\varphi(b)+2.
\end{align*}
The Schwarz lemma and the maximum principle imply that there exists a $w\in\DD$ with $|w|=b$ and
$$
\frac{|f(w)-f(0)|}{|w|}\geq\frac{|f(w_0)-f(0)|}{|w_0|}.
$$
This means that
\begin{align*}
|f(w)|\geq&|f(w)-f(0)|-|f(0)|=|f(0)|+|f(w)-f(0)|-2|f(0)|\\
\geq&\psi(b)-2|f(0)|>\varphi(b)+2-2|f(0)|\geq\varphi(b)=\varphi(|w|).
\end{align*}
This contradicts the maximality of $|f(\zeta)|/\varphi(|\zeta|)$ at $w_0$.
\end{proof}

\begin{lemma}\label{neg2}Assume that $p\in(0,1)$, $q>0$. Then there exist $b\in(0,1)$, $c\in(0,b)$ as
required in Lemma~\ref{neg1}.
\end{lemma}
\begin{proof}We have
$$
\varphi'(t)=\frac{pq}{t^{1-p}(1-t^p)^{q+1}}
$$
which is decreasing for small $t$. To prove the existence of $c$ observe that
$$
\lim_{d\too0^+}\frac{\varphi(d)-\varphi(0)}d=+\infty.
$$
\end{proof}

Observe that $d=k$ for $|z_\ell|$ small enough, $\ell=1,\dots,k$. Let
$h\in\CO(\EE')$ and consider $f(\zeta):=h(\zeta,0,\dots,0)$, $\zeta\in\DD$.
Put $p:=2p_{k+1}$, $q=q_d$. Let
$b,c$ be as in Lemma~\ref{neg2}. It follows from Lemma~\ref{neg1} that
$0\neq|w_0|<c$ cannot be a maximum
of $\DD\ni\zeta\longmapsto|f(\zeta)|/\varphi(|\zeta|)$. In particular $z$
cannot be a maximum of
$\EE'\ni\zeta\longmapsto|h(\zeta)|r_d^{q_d}(\zeta)$ for $0\neq|z_{k+1}|<c$.
\end{proof}

\begin{proof}[Proof of Proposition~\ref{d}]Consider the following lemma.

\begin{lemma}\label{concave}Let $a,c>0$, $t_0\in(0,1)$ be such that
$c\geq1$ or $4a+2c>1+c^2$, $t_0^c>\tau:=a/(a+c)$. Then there exist $b>0$ and
$r\geq1$ such that
$$
h:[0,1]\ni t\longmapsto\frac{t^a}{(r-t)^b}(1-t^c)\in[0,1]
$$
admits its maximum at $t_0$.
\end{lemma}

\begin{proof}Comparing $\varphi:=th'(t)/h(t)$ to zero we get the
equation:
$$
\varphi(r,t)=a+\frac{bt}{r-t}-\frac{ct^c}{1-t^c}=0.
$$
This gives us a formula for $r$:
$$
r(t)=\frac{t(b-bt^c-a+at^c+t^cc)}{-a+at^c+t^cc}.
$$
Observe that:
\begin{align*}
\lim_{t\too\tau^{1/c}}r(t)&=+\infty,\\
\lim_{t\too1}r(t)&=1.
\end{align*}
In order to prove that $r(t)>1$ and that $\varphi(r,\cdot)$ has only one zero
it suffices to show that $r'(t)<0$. We have:
$$
r'(t)=\frac{(-(a+c)t^{2c}+(2a+c-c^2)t^c-a)b
+(at^c+t^cc-a)^2}
{(at^c+t^cc-a)^2}.
$$
It remains to show that the coefficient $\alpha(t^c)$ next to $b$ is negative.
We have:
\begin{align*}
\alpha(\tau)&=\frac{-ac^2}{a+c}<0,\\
\alpha(1)&=-c^2<0,\\
\alpha'(u)&=-2(a+c)u+2a+c-c^2.
\end{align*}
Let $u_0$ be the zero of $\alpha'(u)$. For $c\geq1$ we have $u_0\leq\tau$ and
we
are done. Otherwise $u_0\in(\tau,1)$ and
$4(a+c)\alpha(u_0)=c^2(1+c^2-4a-2c)<0$.
\end{proof}

We may assume that $z_0>0$. Put $a=2p_1$, $c=2p_2$, $t_0=z_0$. Let
$r$ be as in Lemma~\ref{concave}. Putting
$$
h(\zeta)=\frac\zeta{(r-\zeta)^{b/a}}.
$$
completes the proof.
\end{proof}

\begin{remark}\label{partial}The calculations in the proof of Proposition~\ref{d} can be performed using
alternative function families e.g.
\begin{align*}
h_j(\zeta):=&(\zeta+r)^k,\quad r\geq0,\;k=0,1,2,\dots,\\
h_j(\zeta):=&(\zeta+1)^\alpha,\quad \alpha>0,\\
h_j(\zeta):=&(r-\zeta)^\alpha\zeta^k,\quad\alpha<0,r\geq1,k=0,1,2,\dots\\
h_j(\zeta):=&\Big(\frac{\zeta+\delta}{1+\delta\zeta}\Big)\zeta^k,\quad\delta\in[0,1],\;k=0,1,2,\dots.
\end{align*}
However, the author was unable to solve the general case using any of them.
\end{remark}

\noindent{\bf Acknowledgement.} The author is indebted to professor Peter 
Pflug for his valuable comments on the paper.

\end{document}